\documentclass{amsart}

\usepackage{mathrsfs}

\newtheorem{thm}{Theorem}[section]

\newtheorem{prop}[thm]{Proposition}

\theoremstyle{definition}

\theoremstyle{remark}
\newtheorem{oss}[thm]{Remark}

\numberwithin{equation}{section}

\def\R{\mathbb{R}}
\def\P{\mathbb{P}}

\def\wt{\widetilde}

\def\S{\mathbb{S}}
\def\T{\mathbb{T}}

\def\Z{\mathbb{Z}}

\def\EE{\mathcal{E} }

\def\MM{\mathcal{M} }
\def\D{\mathcal{D}  }

\def\LL{\mathcal{L} }
\def\A{\mathcal{A}}

\def\to{\longrightarrow}
\def\sto{\rightarrow}

\def\1{\mathbbm{1}}

\DeclareMathOperator{\Cov}{Cov}

\makeatletter
\newcommand{\vast}{\bBigg@{3.5}}
\newcommand{\Vast}{\bBigg@{5}}
\makeatother
\def\paref#1{(\ref{#1})}

\renewcommand{\(}{\left(}
\renewcommand{\)}{\right)}

\renewcommand{\{}{\left\lbrace}
\renewcommand{\}}{\right\rbrace}

\newcommand{\norm}[1]{\left\lVert#1\right\rVert}

\allowdisplaybreaks

\begin{document}

\title[Random Lipschitz-Killing curvatures]{Random Lipschitz-Killing curvatures: Reduction Principles, Integration by Parts and Wiener chaos}

\author{Anna Vidotto}
\address{Dipartimento di Economia, Universit\`a degli Studi ``G. D'Annunzio'' Chieti-Pescara}
\email{anna.vidotto@unich.it}
\thanks{The author would like to thank Domenico Marinucci for the suggestion of trying to reprove the reduction principle of the Euler-Poincar\'e characteristic via integration by parts. She also would like to thank Maurizia Rossi for many useful comments and  discussions.}

\subjclass[2020]{Primary 60G60,60D05; Secondary 35J05, 60G10, 60G15}

\date{\today}


\keywords{Lipschitz-Killing  curvatures, random eigenfunctions, Wiener chaos expansion, reduction principles}

\begin{abstract}
In this survey we collect some recent results regarding the Lipschitz-Killing curvatures (LKCs) of the excursion sets of random eigenfunctions on the two-dimensional standard flat torus (arithmetic random waves) and on the two-dimensional unit sphere (random spherical harmonics). In particular, the aim of the present survey is to highlight the key role of integration by parts formulae in order to have an extremely neat expression for the random LKCs. 
Indeed, the main tool to study local geometric functionals of random waves on manifold is to exploit their Wiener chaos decomposition and show that (often), in the so-called high-energy limit, a single chaotic component dominates their behavior. Moreover, reduction principles show that  the dominant Wiener chaotic component  of LKCs of random waves' excursion sets at threshold level $u\ne0$ is proportional to the integral of $H_2(f)$, $f$ being the random field of interest and $H_2$ the second Hermite polynomial. This will be shown via integration by parts formulae.
\end{abstract}

\maketitle

\section{Introduction}\label{intro}

The aim of the present survey is to sum up several results presented recently in a number of articles about the local geometry of so-called \emph{random waves}, that are random eigenfunctions of the Laplacian on a compact smooth Riemannian manifold, and giving the reader some evidences for broadly understanding a specific part of this stream of literature. By this, we mean to give some insights of the methodologies and techniques behind most of the proofs, creating, as much as possible, a common thread that unifies the examined works. Lastly, in order to achieve the above-mentioned goal, we will present some alternative proofs of some known results. 

\subsection{Background and notation}\label{B&N}

In order to introduce our framework, let us  fix some preliminary notation. 
Consider a smooth Riemannian manifold $(\MM,g)$ and a random eigenfunction $f_n:\MM\sto \R$ of the Laplacian defined with respect to the Riemannian metric $g$, denoted $\Delta_g$, that is $f_n$ almost surely solves the Helmholtz equation $\Delta_g f_n+\lambda_nf_n=0$, where $-\lambda_n\le0$ is its eigenvalue.
Now fix a level $u\in\R$, we are interested in the geometric properties of the excursion sets of $f_n$, i.e.
$$
\EE_u(f_n,\MM):=\{x\in\MM: f_n(x)\ge u\}\,,
$$
in the high-energy (or high-frequency) limit $\lambda_n\sto\infty$.
Indeed, starting from the seminal work by \cite{Be:77}, the subject has recently attracted great interest, in particular as a consequence of the author conjecturing that, as $\lambda_n\sto\infty$, local geometric functionals of a planar random eigenfunction $f_n$ reflect the behavior of a \emph{typical} deterministic Laplace eigenfunction on any \emph{generic} manifold. 

Indeed, the geometry of excursion sets for \emph{deterministic}, and hence more challenging, eigenfunctions of the Laplacian on a smooth, compact Riemannian manifold have been studied intensively for some time, see among others \cite{Ch:76,Br:78,Ya:82,DF:88} and the recent remarkable articles solving the Yau's conjecture  \cite{Lo:18,LMNN:21}; as a consequence, the \emph{Berry random wave model} attracted several researchers and many articles were written as a consequence of Berry's conjecture, see e.g. \cite{ORW:08,RW:08,Wi:10}.
Another interesting fact that arose from a  work by Berry \cite{Be:02}, is that the fluctuations of the boundary length of $\EE_u(f_n,\D)$ at $u=0$, $\D$ being a smooth subset of  $\R^2$, that is of the so-called \emph{nodal length}, have an unexpected logarithmic order; this was due to a cancellation into the computations whose meaning at that time was considered \emph{obscure} by the author.
This curious phenomenon attracted many researchers to work on local geometric functionals associated to the Berry random wave model on compact manifolds, in particular on their nodal length, see e.g. \cite{KKW:13,Wi:10}, and the cancellation phenomenon was connected to an exact simplification of several terms appearing in the Kac-Rice formula, which is used to compute the so-called \emph{two-point correlation function} of the random geometric functional and hence its variance. Later on, \cite{MW:11, MW:11b} related the cancellation phenomena also to the disappearance, for $u=0$, of the quadratic term in the Hermite expansion of the area of $\EE_u(f_n,\S^2)$, $\S^2$ be the two-dimensional unit sphere, but only in \cite{MPRW:16} the authors elected the Wiener chaos decomposition as fundamental to study the second order behavior of nodal lines on the torus and in general of these local geometric functionals, see Section \ref{LKC and chaos}.

\bigskip

In this survey we want to focus on two cases, when $\MM$ is either the two-dimensional standard flat torus or the two-dimensional unit sphere, denoted by $\T^2$ and $\S^2$ respectively. 
In two dimensions, the so-called Lipschitz-Killing curvatures (see \cite{AT:07}), in the sequel often abbreviated as LKCs, of the excursion sets of the random field $f_n$ characterize its local geometry, those are the Euler-Poincar\'e characteristic $\LL_0^{f_n}(\EE_u(f_n,\MM))$, the boundary length $\LL_1^{f_n}(\EE_u(f_n,\MM))$ and the area $\LL_2^{f_n}(\EE_u(f_n,\MM))$.

\medskip

The above mentioned cancellation phenomena depend on the threshold $u$ and happen at $u=0$ for $\LL_1^{f_n}(\EE_u(f_n,\MM))$ and $\LL_2^{f_k}(\EE_u(f_n,\MM))$, at $u=-1,0,1$ for $\LL_0^{f_n}(\EE_u(f_n,\MM))$.
As a consequence, there is a fundamental difference between what happen at $u=0$, that is for the geometry of the so-called \emph{nodal} sets, and what happen at levels $u\ne0$. In the sequel, we will show explicitly how these cancellations become evident if one consider the Wiener chaos decompositions of the Lipschitz-Killing curvatures, see Remark 3.5 as well as Section \ref{reduction}.

\medskip

Let us briefly give the bibliographic references of the results obtained for LKCs in the case of both random spherical harmonics ($\MM=\S^2$) and arithmetic random waves ($\MM=\T^2$) in the past years.

\subsection{Nodal case}

Marinucci and Wigman \cite{MW:11} studied the variance of the so-called \emph{defect}\footnote{Roughly speaking, the defect of a random eigenfunction is the difference between its positive and negative area.} in the case of $\MM=\S^2$, which is closely  related to the area of the excursion sets at level $u=0$, $\LL_2^{f_k}(\EE_0(f_n,\S^2))$, while in \cite{Ro:15b,Ro:19b} Rossi obtained quantitative central limit theorems for the defect in the case $\MM=\S^d$, $d\ge2$, $\S^d$ being the $d$-dimensional unit-sphere, and showed that its high-energy limit behavior only depend on the odd Wiener chaoses, as the even chaoses vanish at $u=0$. The defect on the two-dimensional standard flat torus was only very recently studied in \cite{KWY:20}, where the authors proved results on the variance and on its spatial distribution.

\medskip

Now let us focus on past works involving $\LL_1^{f_k}(\EE_0(f_n,\MM))$ and related functionals. Number theorists, namely the authors of  \cite{ORW:08} and  \cite{RW:08}, were the first researchers writing on the zeroes of arithmetic random waves. This fact is a direct consequence of the structure of Laplace eigenspaces on the torus, which is inextricably linked to arithmetic considerations, like e.g. that of enumerating lattice points on circles, see Section \ref{ARW}. In particular, the authors of 
\cite{ORW:08} and \cite{RW:08} studied the expectation and variance of the Leray measure of the nodal sets and of the nodal volume, respectively. 
While Rudnick and Wigman \cite{RW:08} only provided a bound for the variance of the volume, the celebrated article \cite{KKW:13} by Krishnapur,  Kurlberg and Wigman gives an exact asymptotic, showing the non-universality of the limit. Such non-universality was reconfirmed by Marinucci, Peccati,  Rossi and Wigman  \cite{MPRW:16}, who provided a non-central limit theorem for the nodal length (see also \cite{DNPR:19}, where the authors showed that the non-universality is preserved by the so-called \emph{phase-singularities}).
On the unit sphere, some preliminary results were obtained already in  \cite{Be:85} and  \cite{Ze:09}, where the expectations of the nodal lengths of the long and small energy window random functions were computed. However, the first results on local geometric functionals on $\S^2$ that can be compared to the ones obtained in \cite{KKW:13}  on $\T^2$ were presented by Wigman \cite{Wi:10}, who computed the expectation and variance of the nodal length of the spherical harmonics. Moreover, Marinucci, Rossi and Wigman \cite{MRW:20} obtained a central limit theorem, in the high frequency limit (see also the survey \cite{Ro:15b} as well as the interesting monograph  \cite{MP:11}). 

\medskip

The study of the Euler-Poincar\'e characteristic at level $u=0$ is still open; what  is known is that a cancellation is occurring at levels $u=-1,0,1$, see \cite{CM:18,CMR:20} for the results on the sphere and on the torus respectively.

\subsection{Non-zero levels}

Let us now consider \emph{non-zero level} Lipschitz-Killing curvatures. 
Regarding the excursion area, that is $\LL_2^{f_n}(\EE_u(f_n,\MM))$ when $u\ne0$. In the case of the two dimensional sphere, Marinucci and Wigman \cite{MW:11b} computed the variance and obtained a CLT, while Marinucci and Rossi \cite{MR:15} extended the results to $\S^d$, $d\ge2$. For arithmetic random waves, analogous results can be found in the recent article \cite{CMR:20}.

\medskip

Let us now focus on the boundary length for $u\ne0$: the computation of the variance and a CLT on $\S^2$ can be found once again in Rossi PhD thesis \cite{Ro:15b}, while analogous results on $\T^2$ can be found in [Remark 2.4]\cite{MPRW:16}, as well as \cite{Ro:19} and \cite{CMR:20}, showing that the universality is preserved in the non-nodal case. Indeed, it is important to underline that, while for the spherical harmonics both the nodal length and the non-zero-level length converge in distribution to a normal random variable, on the torus the nodal length converges in law to a linear combination of independent chi-squares (where  the coefficients depend on the  chosen subsequence of energy levels) and the non-zero-level length to a Gaussian random variable.

\medskip

Regarding the EPC at $u\ne0$, complete results, namely on the variance and on the high-energy limit distribution, with a quantitative CLT, for spherical harmonics can be found in \cite{CM:18}, while for arithmetic random waves once again in \cite{CMR:20}.

Related results on the same type of functionals but on other manifolds can be found in \cite{RWY:16},  \cite{Ca:19},  \cite{Ma:17},  \cite{Ma:17b},  \cite{BM:17},  \cite{PR:18},  \cite{RW:18},  \cite{Ma:18}, \cite{MP:10},  \cite{MW:14},  \cite{CMR:18},  \cite{MR:15},  \cite{CMW:16},  \cite{CMW:16b},  \cite{BDMP:16},   \cite{CW:17},  \cite{To:19},  \cite{To:20}, \cite{NPR:19,PV:20,Vi:21}, \cite{BCW:17},\cite{No:20}.

\subsection{Motivations and aim of the survey}
 
Even if the Wiener chaos expansion was already used to study random geometric functionals, for instance in \cite{KL:10, AL:13, EL:16}, the work by Marinucci, Peccati, Rossi and Wigman \cite{MPRW:16} was fundamental in presenting the remarkable idea that has allowed to prove limit theorems for nodal (and then also non-nodal) statistics of random waves on compact manifolds: first one has to derive their chaos decompositions, then prove that, most of the time, a single chaotic projection is dominant. On the sphere and on the torus, the dominating chaotic projections can often be represented as explicit functionals of a finite collection of independent Gaussian coefficients and one can make use of the standard CLT or its possible generalizations\footnote{Another well-known strategy could be the use  of the so-called Fourth Moment Theorem by \cite{NP:05} and in its quantitative form by \cite{NP:09}. This can be crucial in situations in which  the standard CLT cannot  be applied.}. Moreover, first in \cite{CM:18} and then in \cite{MRW:20, CMR:20}, the authors were able to prove some reduction theorems in which they show that the dominant chaotic component could even have a neater  expression, given by a deterministic function of the threshold $u$ times the integral of a single Hermite polynomial evaluated on the field $f_n$. Some immediate reduction principles were already proved via integration by parts formulae, in particular the Green identity, for the boundary length at $u\ne0$ in \cite{Ro:15b, Ro:19}, in order to prove CLTs via usual tools and show the exact cancellation of the second chaos at $u=0$. However, the extension to the Euler-Poincar\'e characteristics in \cite{CM:18} and \cite{CMR:20} rely on a different approach requiring more involved computations. 
Regarding the reduction principle for the nodal length on the sphere in \cite{MRW:20} (see also \cite{Vi:21} for analogous results on the plane), the situation is not directly comparable since it is an asymptotic full correlation result, i.e. the authors were able to prove that the behavior of the single dominant chaotic component is asymptotically the same of a simpler statistic, the so-called \emph{sample trispectrum}.

\medskip

The aim of the present survey is to present the fundamental steps that, after the seminal work  \cite{MPRW:16}, allowed to prove  high-energy limit results in the above mentioned  works, focussing on Wiener chaos expansions and reduction principles. 
In particular, we will not  present the asymptotic  results  regarding the variance and   the  distributional limit of the three LKCs, but only show the crucial role of the Wiener chaos expansion  of these local geometric functionals, together with their reduction theorems that can be beautifully obtained via integration by parts formulae for the non-zero level case,  $u\ne0$.
Finally, we will show that the reduction principles obtained in \cite{CM:18} and \cite{CMR:20}, that is for the Euler-Poincar\'e characteristics on the sphere and on the torus respectively, can also be proved via integration by parts formulae.

\begin{oss}
In this survey, by reduction principles we mean not only the fact that a single chaotic component dominates the functional of interest, but also that this single chaotic component can be written as a function of the level $u\ne0$ times the integral of $H_2(f_n)$.
Moreover, it is important to stress that the reduction principles for the boundary length on both manifolds and for the Euler-Poincar\'e characteristic on the torus hold for each $n$, namely are non-asymptotic results. The high-energy asymptotic regime starts to play a fundamental role when one wants to prove that a single component in the Wiener chaos expansion dominates the Lipschitz-Killing curvature of interest. See also Remarks \ref{rem-NA-2} and \ref{rem-NA-3}.
\end{oss}

\begin{oss}
Note that in \cite{No:20}, an integration by parts formula for some general functionals of independent random field is presented. However, this formula does not cover the case of the Euler-Poincar\'e characteristic.
\end{oss}

\begin{oss}
At the end of this introduction, it is worth pointing out in which sense the geometric functionals considered here have a \emph{local} nature, as often mentioned before: the crucial point is that the three Lipschitz-Killing curvatures satisfy some additivity properties. Indeed, it is always possible to exploit the fact that, if $A,B$  are two closed convex subsets of $\MM$ and $A\cap B=\emptyset$, then $\LL_j(A\cup B)=\LL_j(A)+\LL_j(B)$, $j=0,1,2$. 
This property fails for so-called \emph{global} geometric quantities associated with excursion sets, whose study becomes more challenging. To mention some results of this type, Nazarov and Sodin \cite{NS:16} studied the expectation of the number of connected components of nodal sets of generic Gaussian random functions of several real variables, while, for the so-called random band-limited functions, in \cite{SW:16} is showed that topologies and nestings of the zero and nodal sets have universal laws of distribution.
\end{oss}

\paragraph{\bf Plan of the survey} In Section 2 we briefly present the construction of arithmetic random waves and random spherical harmonics, and we  set  most of the notations. In Section 3 we give a brief compendium on Wiener chaos and we show the chaotic decompositions of the three LKCs. Finally, in Section 4 we present the reduction principles together with their proofs, in particular giving an alternative proof for the reduction formula of the Euler-Poincar\'e characteristic.

\section{Random eigenfunctions}

In Section \ref{B&N} we gave a brief definition of random eigenfunctions of the Laplacian. Since we will focus on $\MM$ being either the two-dimensional standard flat torus or the two-dimensional unit sphere, here we want  to introduce in more detail random eigenfunctions of the Laplacian on these two specific smooth compact Riemannian manifolds or, more precisely, we will present how these eigenfuctions are constructed in the Gaussian random framework. 
In particular,  we will see that, in the high-energy limit, $\text{i.e.}$ when the  eigenvalues diverge, the covariance  structure of the Gaussian random eigenfunction converges in some sense to the one of the Berry random wave model on the Euclidean plane, whose covariance kernel is, for $x,y\in \R^2$,
\begin{equation}\label{cov-berry}
\Cov\(f(x),f(y)\)=J_0\(\norm{x-y}\),
\end{equation}
$J_0$  being the Bessel function of  order $0$ (see \cite[Section 1.71]{Sz:75})  and $\norm{\cdot}$ denoting the Euclidean norm. It is worth stressing that this model, according to Berry, predict the local behavior of deterministic eigenfunction on generic chaotic surfaces for large eigenvalues, see \cite{Be:77, Be:02}.

\subsection{Arithmetic random waves}\label{ARW}

It is a standard fact that the Laplacian eigenvalues for the two-dimensional standard flat torus $\T^2:=\R^2/\Z^2$ are of the form $-\lambda_n=-4\pi^2n$, where $n$  is an integer that can be written as a sum of two squares, $\text{i.e.}$
$$
n\in S := \{ n = a^2 + b^2 : a,b\in \Z \}\,,
$$ 
this is why the eigenspaces of the Laplacian on $\T^2$ are related to the theory of lattice points on circles.
In order to introduce the eigenspace associated to each $\lambda_n$, we need to  define the set of frequencies  
\begin{equation}\label{Lambda}
\Lambda_n := \lbrace \xi\in \mathbb Z^2 : \| \xi\|= \sqrt n\rbrace\,, \qquad n \in S\,,
\end{equation}
denoting $\mathcal N_n$ its cardinality, so that $\mathcal N_n$ is the multiplicity of $\lambda_n$. 
In particular, via the set $\Lambda_n$ in \paref{Lambda} and denoting $\delta_z$ the Dirac mass at $z\in \R^2$, one can define a probability measure $\mu_n$ on the unit circle $\mathbb S^1\subset \R^2$
$$
\mu_n := \frac{1}{\mathcal N_n} \sum_{\xi\in \Lambda_n} \delta_{\xi/\sqrt n}\,,
$$ 
see \cite{KKW:13} for a more detailed discussion.

The eigenspace $\mathcal E_n$ associated with $\lambda_n$ is spanned  by the $L^2$-orthonormal set of functions $\{e^{i2\pi\langle\xi,\cdot\rangle}\}_{\xi\in\Lambda_n}$, so that, for $n\in S$, the arithmetic random wave of order $n$ is a Gaussian randomization of functions living in the eigenspace $\mathcal E_n$, or, being more precise, a random linear combination of the following form, see \cite{RW:08}:
\begin{equation}\label{defARW}
f_n(x) := \frac{1}{\sqrt{\mathcal N_n}} \sum_{\xi\in \Lambda_n}a_\xi e^{i2\pi \langle \xi, x\rangle},\qquad x\in \T^2,
\end{equation}
where $\lbrace a_\xi \rbrace_{\xi\in \Lambda_n}$ is a family of identically distributed standard complex Gaussian random variables, defined on some probability space $(\Omega, \mathcal F, \mathbb P)$, and independent except for the relation $\overline{a_\xi} = a_{-\xi}$ that ensures $f_n$ to be real.  
Equivalently, the random field $f_n$ can be defined as the centered Gaussian function on $\T^2$ whose covariance kernel  is, for $x,y\in \T^2$,
\begin{equation}\label{cov_int}
\Cov(T_n(x), T_n(y)) = \frac{1}{\mathcal N_n} \sum_{\xi\in \Lambda_n} e^{i 2\pi\langle \xi, x-y\rangle } 
\end{equation}
and by this alternative definition it is clear that the law of $f_n$ is invariant under translations, namely that $f_n$ is stationary.
Note there exists a density one subsequence $\lbrace n_j\rbrace_j\subset S$ of energy levels such that, as $j\sto +\infty$, 
$$
\mu_{n_j} \Rightarrow d\theta/2\pi,
$$
$d\theta$ denoting the uniform measure on $\mathbb S^1$, see \cite{FKW:06}. 
From \paref{cov_int} we have, for $x,y\in \T^2$, as $j\sto +\infty$, 
\begin{equation*}
\Cov(T_{n_j}(x/2\pi \sqrt{n_j}), T_{n_j}(y/2\pi\sqrt{n_j})) = \int_{\mathbb S^1} e^{i \langle \theta, x-y\rangle} d\mu_{n_j}(\theta) \to J_0(\|x-y\|),
\end{equation*}
$J_0$ still denoting the Bessel function of order zero, showing the convergence to the Berry Random Wave Model \eqref{cov-berry}.
A partial classification of the others possible weak-$\star$ limits of subsequences of $\lbrace \mu_n\rbrace_{n\in S}$ can be found in \cite{KW:17}. 

\begin{oss}\label{rem-NA-2}
The result for $\MM=\T^2$ presented in this  survey are non-asymptotic,  however, they are the preamble of most of the high-energy limit results obtained in the works mentioned in the introduction. On the torus, by high-energy limit one means the asymptotic behavior as $n\sto\infty$ such that $\mathcal N_n\sto\infty$, which is not granted as for the spherical case (see the next section). Anyway we  can say that  $\mathcal N_n$  grows   on average as $\sqrt {\log n}$, see \cite{La:08}. As a consequence, talking about high-energy limit on the torus, one has to take the extra assumption that $\mathcal N_n$ grows to infinity with $n$. This assumption will be tacit in the sequel. See also Remark \ref{rem-NA-3}.
\end{oss}

\subsection{Random spherical harmonics}\label{RSH}
 
On the two-dimensional unit sphere $\mathbb S^2$, the Laplacian eigenvalues are of the form $-\lambda_n=-n(n+1)$, where $n\in \mathbb N$, the multiplicity of the $n$-th eigenvalue is $2n+1$ and the family $\lbrace Y_{n,m}\rbrace_{m=-n,\dots , n}$ of deterministic functions, which are a base  of the so-called spherical harmonics \cite[Section 3.4]{MP:11}, represents an orthonormal basis for the eigenspace $\mathcal E_n$ corresponding to $\lambda_n$.
As a consequence, one can choose an arbitrary $L^2$-orthonormal basis for $\mathcal E_n$ and construct the $n$-th random spherical harmonic on $\S^2$ in the following way, see \cite{Wi:10}:
\begin{equation}\label{defRSH}
f_n(x) := \sqrt{\frac{4\pi}{2n+1}} \sum_{m=-n}^n a_{n,m} Y_{n,m}(x),\qquad x\in \mathbb S^2,
\end{equation}
where $\lbrace a_{n,m} \rbrace_{m=-n, \dots, n}$ is a family of identically distributed standard complex Gaussian random variables, and independent except for the relation $\overline{a_{n,m}} = (-1)^n a_{n,-m}$ that ensures $f_n$ to be real. 
Also on the sphere, the random field $f_n$ can be equivalently defined as the centered Gaussian function on $\S^2$ whose covariance kernel is, for $x,y\in \mathbb S^2$,
$$
\Cov(f_n(x)), f_n(y)) = P_n(\cos d(x,y)),
$$
where $P_n$ denotes the $n$-th Legendre polynomial, see \cite[Section 13.1.2]{MP:11}, and $d(x,y)$ the geodesic distance between the two points $x,y\in\S^2$. 
Thanks to Hilb's asymptotic formula \cite[Theorem 8.21.12]{Sz:75}, which states that, \emph{uniformly} for  $\theta\in [0, \pi - \varepsilon]$ ($\varepsilon >0$), as $n\sto +\infty$, 
\begin{equation}\label{Hilb}
P_n(\cos \theta) \sim \sqrt{\frac{\theta}{\sin \theta}} J_0((n+1/2)\theta),
\end{equation}
we have again some heuristics regarding the above-mentioned Berry's conjecture \eqref{cov-berry}.

\begin{oss}
The reader probably already noticed the unified notation $f_n, -\lambda_n$ for the Gaussian random eigenfunction and its eigenvalue, respectively, that lives either on the two dimensional unit sphere or on the two dimensional standard flat torus. Since the results presented in this survey are independent of the manifold, $f_n$ will denote, from now on, both the random fields. In particular, if $\MM=\T^2$ one has to think about arithmetic random waves and the fact that $n\in S$ will be implied, otherwise, if $\MM=\S^2$ one has to think about random spherical harmonics.
\end{oss}

\section{Lipschitz-Killing curvatures and Wiener chaos}\label{LKC and chaos}

Let us recall once again  that for a random field with a two-dimensional domain, Lipschitz-Killing curvatures are the three quantities that characterize any local geometric functional associated with its excursion sets, for a detailed discussion see \cite[Ch.6, Section 6.3]{AT:07} and in particular Theorem  6.3.1.
Moreover, several recent articles have shown that the Lipschitz-Killing curvatures for the excursion sets of both random spherical harmonics and arithmetic random waves are often dominated, in the high-energy limit $n\sto\infty$, by a single Wiener chaotic component. 
More precisely, for the Lipschitz-Killing curvatures of the excursion sets at level $u\ne0$, this dominant chaotic component is the projection onto the second  order  Wiener chaos and it can be written, thanks to the reduction principles proved in \cite{Ro:15b, CM:18, CMR:20}, as a simple explicit function of the threshold parameter $u$ times the centered norm of these random fields, see also Section \ref{reduction}; this is why its disappearance results in a smaller order variance and on a different limiting behavior. For this reason, in Section \ref{wc} we will present a short compendium on Wiener chaos and then in Section \ref{ce-LPCs} we will show results regarding the Wiener chaotic decompositions for the three Lipschitz-Killing curvatures, namely the excursion area, the boundary length and the Euler Poincar\'e characteristic, that were obtained by the various authors that worked on this topic.

We stress that, regarding \emph{nodal} Lipschitz-Killing curvatures, the dominant term is not always a single chaotic component and it can be of different order, for instance the fourth chaos for the nodal length and odd chaoses for the defect.

\subsection{Wiener chaos}\label{wc}

As previously remarked, since  \cite{MPRW:16}, the Wiener chaos decomposition of local geometric functionals such as LKCs plays a fundamental role in understanding their asymptotic behavior, so here we introduce the concept of Wiener chaos both on the torus and on the sphere, trying as much as possible to unify the framework. For a complete discussion on Wiener chaos see  [Section 2.2]\cite{NP:12} and the references therein. From now on, $\MM$ will denote either the two-dimensional sphere $\S^2$ or the two-dimensional torus $\T^2$.

Denote by $\{H_k\}_{k\ge 0}$ the sequence of Hermite polynomials on $\mathbb{R}$; these polynomials are defined recursively as follows: $H_0 \equiv 1$ and
 $$H_{k}(t) = tH_{k-1}(t) - H'_{k-1}(t), \qquad k\ge 1.$$
Recall that $\mathbb{H} := \{(k!)^{-1/2} H_k, k\ge 0\}$ forms a complete orthonormal system in the space of square integrable real functions $L^2(\gamma)$ with respect to the standard Gaussian density $\gamma$ on the real line, see [Section 1.4]\cite{NP:12}.

Arithmetic Random Waves \eqref{defARW} are generated from a family of complex-valued Gaussian random variables $\{a_\xi\}_{\xi\in \mathbb{Z}^2}$, defined on  $(\Omega, \mathscr{F}, \mathbb{P})$ and verifying the following properties: {\bf (1)} every $a_\xi$ has the form $a_\xi= \Re(a_\xi) + i \Im(a_\xi)$ where $\Re(a_\xi)$ and $\Im(a_\xi)$ are two independent real-valued, centred, Gaussian random variables with variance 1/2,
{\bf (2)} the $a_\xi$'s are stochastically independent, save for the relations $a_{-\xi} = \overline{a_\xi}$ in particular making $f_n$ real-valued.
In the case $\MM=\T^2$, let us define the space ${\mathcal A}$ to be the closure in $L^2(\mathbb{P})$ of all real finite linear combinations of random variables $\zeta$ of the form $$\zeta = z \, a_\xi + \overline{z} \, a_{-\xi}\,,$$ where $\xi\in \mathbb{Z}^2$ and $z\in \mathbb{C}$.

Random spherical harmonics \eqref{defRSH} are generated from a family of complex-valued Gaussian random variables $\{a_{\ell,m} : \ell = 0, 1, 2, \dots, m=-\ell, \dots, \ell\}$ such that {\bf (a)} every $a_{\ell,m}$ has the form $x_{\ell,m}+iy_{\ell,m}$, where $x_{\ell,m}$ and $y_{\ell,m}$ are two independent real-valued Gaussian random variables with mean zero and variance $1/2$; {\bf (b)} $a_{\ell,m}$ and $a_{\ell',m'}$ are  independent whenever $\ell \ne \ell'$ or $m'\notin \{m,-m\}$, and {\bf (c)} $\overline{a_{\ell,m}} =(-1)^\ell a_{\ell, -m}$. In the case $\MM=\S^2$, define the space ${\A}$ to be the closure in $L^2(\mathbb{P})$ of all real finite linear combinations of random variables $\zeta$ of the form $$\zeta = z \, a_{\ell,m} + \overline{z} \, (-1)^\ell a_{\ell, -m},\qquad z\in \mathbb{C}\,.$$ 

In both cases, $\A$ is a real centered Gaussian Hilbert subspace
of $L^2(\mathbb{P})$.

{\rm Let us fix now an integer $q\ge 0$; the $q$-th Wiener chaos $C_q$ associated with $\A$ is defined as the closure in $L^2(\mathbb{P})$ of all real finite linear combinations of random variables of the type
$$
H_{p_1}(\xi_1)\cdot H_{p_2}(\xi_2)\cdots H_{p_k}(\xi_k)
$$
for $k\ge 1$, where the integers $p_1,...,p_k \geq 0$ satisfy $p_1+\cdots+p_k = q$, and $(\xi_1,...,\xi_k)$ is a standard real Gaussian vector extracted
from $\A$ (in particular, $C_0 = \mathbb{R}$).}

Taking into account the orthonormality and completeness of $\mathbb{H}$ in $L^2(\gamma)$ (see e.g. \cite[Theorem 2.2.4]{NP:12}), it is possible to prove that $C_q \,\bot\, C_m$ in  $L^2(\mathbb{P})$ for every $q\neq m$, and moreover
\begin{equation*}
L^2(\Omega, \sigma(\A), \mathbb{P}) = \bigoplus_{q=0}^\infty C_q,
\end{equation*}
that is, every real-valued functional $F$ of $\A$ can be (uniquely) represented as a series, converging in $L^2$, of the form
\begin{equation}\label{e:chaos2}
F = \sum_{q=0}^\infty F[q],
\end{equation}
where $\text{\rm Proj}[F | q]$ stands for the projection of $F$ onto $C_q$, and the series converges in $L^2(\mathbb{P})$. Plainly, $\text{\rm Proj}[F | 0]= \mathbb E[F]$.

\medskip

In the sequel, for $i,j=1,2$, we will denote by $\partial_{i} f_n(x)=\partial_{i} f_n(x_1,x_2)$ the partial derivative with respect to $x_i$ and denote $\partial_{ij} f_n(x)=\partial_{ij} f_n(x_1,x_2)$ the second partial derivative with respect to $x_i$ and $x_j$. For $\MM=\S^2$, note that $x=(x_1,x_2)=(\theta,\phi)$, $\theta\in(0,\pi), \phi\in[0,2\pi)$, and in this system of coordinates the gradient is given by  $\nabla=\(\frac{\partial}{\partial  \theta}, \frac{1}{\sin\theta}\frac{\partial}{\partial \phi}\)$.
The random fields $f_{n},\partial_{j} f_n,\partial_{ij} f_n$ viewed as collections of
Gaussian random variables
indexed by $x\in \MM$ are all lying in ${\mathcal A}$, i.e. for every $x\in \MM$ we have
\begin{equation*}
f_{n}(x),\, \partial_{j}f_{n}(x), \, \partial_{ij}f_{n}(x) \in \mathcal A.
\end{equation*}

\subsection{Chaotic expansions of Lipschitz-Killing curvatures}\label{ce-LPCs}

We recall that, from now on, $\MM$ denotes either $\T^2$ or $\S^2$, and it will not be specified when not needed.
The Lipschitz-Killing curvatures are finite-variance functionals of $\mathcal A$, hence applying (\ref{e:chaos2}) we get the series expansion 
$$
\LL_k^{f_n}(\EE_u(f_n,\MM))=\sum_{q=0}^\infty \operatorname{Proj}(\LL_k^{f_n}(\EE_u(f_n,\MM))|q). 
$$
Let us be more precise.

\medskip

The excursion area has the following integral representation, see also [Section 3]\cite{MW:11},
\begin{equation}
\LL_2^{f_n}(\EE_u(f_n,\MM)) = \int_{\MM} \mathbf{1}_{\lbrace f_n(x)\ge u\rbrace}\,dx
\end{equation}
which guarantees that $\LL_2^{f_n}(\EE_u(f_n,\MM))\in L^2_{\mathcal A}(\mathbb P)$. In \cite{MW:11} and \cite{CMR:20} one  can find the following result, which is proven in detail in [Section 3]\cite{MW:11}.
\begin{prop}[\cite{MW:11,CMR:20}]
For every $n$ such that $f_n$ is an eigenfunction of $\Delta_g$ with eigenvalue $\lambda_n$ and $u\in \mathbb R$, the chaotic decomposition of $\LL_2^{f_n}(\EE_u(f_n,\MM))$ is given by
\begin{equation}
\LL_2^{f_n}(\EE_u(f_n,\MM)) = \sum_{q=0}^{+\infty} \frac{\gamma_q(u)}{q!} \int_{\MM} H_q(f_n(x))\,dx,
\end{equation}
where $\gamma_q(u):= H_{q-1}(u)\phi(u)$, and the convergence of the series is in $L^2(\mathbb P)$.
\end{prop}

\medskip

The boundary length has the following \emph{formal} integral representation, see [Section 7.2.1]\cite{Ro:15b}
\begin{equation}\label{formalBL}
\LL_1^{f_n}(\EE_u(f_n,\MM)) = \frac12 \int_{\MM} \delta_u(f_n(x)) \| \nabla f_n(x)\|\,dx,
\end{equation}
where $\delta_u$ is the Dirac mass in $u$, and $\nabla f_n$ is the gradient of $f_n$. 
Let us introduce two collections of coefficients
$\{\alpha_{2n,2m} : n,m\geq 1\}$ and $\{\beta_{l}(u) : l\geq 0\}$, that are needed in order to state  the chaotic expansion of $\LL_1^{f_n}(\EE_u(f_n,\MM))$ and  are related to the Hermite expansion of the norm $\| \cdot
\|$ in $\R^2$ and the (formal) Hermite expansion of the Dirac mass $ \delta_u(\cdot)$ respectively, see \cite{Ro:15b}.
These are given by
\begin{equation}\label{e:beta}
\beta_{l}(u):= H_{l}(u)\phi(u),
\end{equation}
where $H_{l}$ still denotes the $l$-th Hermite polynomial and
\begin{equation}\label{e:alpha}
\alpha_{2n,2m}:=\sqrt{\frac{\pi}{2}}\frac{(2n)!(2m)!}{n!
m!}\frac{1}{2^{n+m}} p_{n+m}\left (\frac14 \right),
\end{equation}
where for $N\in \mathbb N$ and $x\in \R$
\begin{equation*}
\displaylines{ p_{N}(x) :=\sum_{j=0}^{N}(-1)^{j}\cdot(-1)^{N}{N
\choose j}\ \ \frac{(2j+1)!}{(j!)^2} x^j, }
\end{equation*}
the ratio $\frac{(2j+1)!}{(j!)^2}$ being the so-called swinging factorial
restricted to odd indices. 

\begin{prop}[\cite{Ro:15b, MPRW:16, CMR:20}]\label{lem-ce1}
For every $n$ such that $f_n$ is an eigenfunction of $\Delta_g$ with eigenvalue $\lambda_n$ and $u\in \mathbb R$ the chaotic expansion of $\LL_1^{f_n}(\EE_u(f_n,\MM))$ is
\begin{eqnarray*}
\LL_1^{f_n}(\EE_u(f_n,\MM)) &=& \frac12 \sqrt{\frac{\lambda_n}{2}}\sum_{q=0}^{+\infty}\sum_{u=0}^{q}\sum_{k=0}^{u}
\frac{\alpha _{2k,2u-2k}\beta _{q-2u}(u)
}{(2k)!(2u-2k)!(q-2u)!}\times\\
\nonumber
&&\times \int_{\mathbb T}H_{q-2u}(f_n(x))
H_{2k}(\wt \partial_1f_{n}(x))H_{2u-2k}(\wt \partial_2f_{n}(x))\,dx,
\end{eqnarray*}
where the convergence of the series is in $L^2(\P)$, and $\wt \partial_j f_{n}$, $j=1,2$,x denotes normalized first derivatives.
\end{prop}

\medskip

The Euler-Poincar\'e characteristic has the following \emph{formal} representation
\begin{equation}\label{formalEPC}
\LL_0^{f_n}(\EE_u(f_n,\MM)) = \int_{\mathbb T} \text{\rm det}(\nabla^2 f_n(x)) \mathbf{1}_{\lbrace f_n(x)\ge u\rbrace} \delta_0(\nabla f_n(x))\,dx,
\end{equation}
where $\nabla^2 f_n$ is the Hessian matrix of $f_n$, and abusing notation $\delta_0$ denotes the Dirac mass in $(0,0)$.

The following result presents the chaotic expansion of $\LL_0^{f_n}(\EE_u(f_n,\MM))$, note that \cite{CM:18} and \cite{CMR:20} do not give explicit expressions for chaotic coefficients but those corresponding to the zero-th and second Wiener chaoses.
\begin{prop}[\cite{CM:18, CMR:20}]\label{prop-chaos-EPC}
For $n\in S$ and $u\in \mathbb R$, the chaotic expansion of $\LL_0^{f_n}(\EE_u(f_n,\MM))$ is
 \begin{equation}\label{exp_L}
  \begin{split}
  \LL_0^{f_n}(\EE_u(f_n,\MM))
    =& 2\lambda_n\sum_{q=0}^{+\infty} \sum_{a+b+c+2d+2e=q} \frac{\eta^{(n)}_{a,b,c}(u)}{a!b!c!}\frac{\beta_{2d} \beta_{2e}}{(2d)!(2e)!} \int_{\mathbb T} H_{a}\left (\frac{\partial_{11}f_n(x)}{\kappa_3} \right ) \\
    &\times H_b \left (\frac{\partial_{12}f_n(x)}{\kappa_4}\right )
     H_c \left (\frac{\partial_{22}f_n(x)}{\kappa_5} - \frac{\kappa_2}{\kappa_5 \kappa_3} \partial_{11}f_n(x) \right )
    H_{2d}\left (\frac{\partial_1 f_n(x)}{\kappa_1}\right )\\
    &\times  H_{2e}\left (\frac{\partial_2 f_n(x)}{\kappa_1} \right )\,dx,
  \end{split}
  \end{equation}
 for some coefficients $\eta^{(n)}_{a,b,c}(u)\in \mathbb R, a,b,c\in \mathbb N$, where the series converges in $L^2(\mathbb P)$,
  \begin{equation}\label{e:beta0}
  \beta_q := \beta_q(0) = \phi(0) H_q(0)
  \end{equation}
  as defined in (\ref{e:beta}),
  and $\kappa_1,\dots,\kappa_5$ are for $\MM=\S^2$,
  \begin{eqnarray*}
&&\kappa_1=\frac{\sqrt{\lambda }_{\ell }}{\sqrt{2}}\qquad \kappa_2=
\frac{\sqrt{\lambda _{\ell }}(\lambda _{\ell }+2)}{2\sqrt{2}\sqrt{3\lambda
_{\ell }-2}}\qquad \kappa_3=\frac{\sqrt{\lambda _{\ell }}\sqrt{3\lambda
_{\ell }-2}}{2\sqrt{2}} \\
&&\kappa_4=\frac{\sqrt{\lambda _{\ell }}\sqrt{\lambda _{\ell }-2}}{2
\sqrt{2}}\qquad \kappa_5=\frac{\lambda _{\ell }\sqrt{{\lambda _{\ell }-2}
}}{\sqrt{3\lambda _{\ell }-2}}
\end{eqnarray*}
  while for $\MM=\T^2$,
  \begin{flalign*}
&\kappa_1=\sqrt{\frac{\lambda_n}{2}} \qquad \kappa_2=\frac{\lambda_n}{2\sqrt2}\frac{1-\hat{\mu}_n(4)}{\sqrt{3+\hat{\mu}_n(4)}} \qquad \kappa_3=\frac{\lambda_n}{2\sqrt2}\sqrt{3+\hat{\mu}_n(4)}\\
& \kappa_4=\frac{\lambda_n}{2\sqrt2}\sqrt{1-\hat{\mu}_n(4)} \qquad \kappa_5=\lambda_n\frac{\sqrt{1+\hat{\mu}_n(4)}}{\sqrt{3+\hat{\mu}_n(4)}}\,,
\end{flalign*}
where $\hat{\mu}_n(4)$ is the fourth Fourier coefficient of $\mu_n$.
\end{prop}

\begin{oss}
The two formal representations \eqref{formalBL} and \eqref{formalEPC} are justified by the use of the  following $\varepsilon$-approximating random variables, $\varepsilon>0$, 
\begin{flalign}
&\mathcal \LL_1^{f_n,\varepsilon}(\EE_u(f_n,\MM)) := \frac12 \int_{\MM} \frac{1}{2\varepsilon}\mathbf{1}_{[u-\varepsilon, u+\varepsilon]}(f_n(x))\| \nabla f_n(x)\|\,dx\label{approx1}\\
&\mathcal \LL_0^{f_n,\varepsilon}(\EE_u(f_n,\MM)) = \int_{\mathbb T}  \text{\rm det}(\nabla^2 f_n(x)) \mathbf{1}_{\lbrace f_n(x)\ge u\rbrace} \frac{1}{(2\varepsilon)^2}\mathbf{1}_{[-\varepsilon, \varepsilon]^2}(\nabla f_n(x))\,dx\label{approx0}
\end{flalign}
that, uniformly in $n$, converge to the first and the zero-th Lipschitz-Killing curvature  respectively. The chaotic expansions are then computed for \eqref{approx1} and \eqref{approx0}, and then extended to the original LKCs, see the detailed \cite[Section 4.2]{CMR:20} and the references therein.
\end{oss}
  

\section{Reduction principles and integration by parts formulae}\label{reduction}

For the case $u\ne0$, the three Lipschitz-Killing Curvatures are dominated by their second-order chaotic components, that are given by
\begin{flalign}
 \operatorname{Proj}[ \LL_0^{f_n}(\EE_u(f_n,\MM)) |2]
    &= 2\lambda_n \sum_{a+b+c+2d+2e=2} \frac{\eta^{(n)}_{a,b,c}(u)}{a!b!c!}\frac{\beta_{2d} \beta_{2e}}{(2d)!(2e)!} \int_{\mathbb T} H_{a}\left (\frac{\partial_{11}f_n(x)}{\kappa_3} \right ) \notag\\
    &\times H_b \left (\frac{\partial_{12}f_n(x)}{\kappa_4}\right )
     H_c \left ( \frac{\partial_{22}f_n(x)}{\kappa_5} - \frac{\kappa_2}{\kappa_5 \kappa_3} \partial_{11}f_n(x)\right )\notag\\
    &\times  H_{2d}\left (\frac{\partial_1 f_n(x)}{\kappa_1}\right ) H_{2e}\left (\frac{\partial_2 f_n(x)}{\kappa_1} \right )\,dx\,,\label{2ndchaos-EPC}\\
 \operatorname{Proj}[ \LL_1^{f_n}(\EE_u(f_n,\MM)) |2]&= \sqrt{\frac{\lambda_n}{2}}\Big (\frac{\beta_2 \alpha_{00}}{2!}\int_{\MM} H_2\(f_n(x)\)\,dx \notag\\
&+ \frac{\beta_0 \alpha_{20}}{2!}\int_{\MM} \{H_2\(\wt\partial_1f_n(x)\)+H_2\(\wt\partial_2f_n(x)\)\} dx\,,\label{2ndchaos-BL}\\
 \operatorname{Proj}[ \LL_2^{f_n}(\EE_u(f_n,\MM)) |2]&= \frac{1}{2} u\, \phi(u)\,\int_{\MM} (f_n(x)^2 - 1)\,dx \,\notag;
\end{flalign}
note that for the boundary length and the Euler-Poincar\'e characteristic, the second chaos seems to depend on the derivatives of the field $f_n$.
However, in \cite{Ro:15b, CM:18}  on the sphere, in \cite{CMR:20} on the torus, it was shown that  all the expressions of the projections onto the second order Wiener chaos can  be reduced to the following beautiful formula,  which involves a deterministic function of the level $u$ and the integral of the second Hermite polynomial $H_2$, evaluated on the field $f_n$:
\begin{flalign}
&\operatorname{Proj}[\LL_k^{f_n}(\EE_u(f_n,\MM))|2] = c_k(u) \left (\sqrt{\frac{\lambda_n}{2}} \right )^{2-k} \int_{\MM} H_2(f_n(x))\,dx +O_{L^2(\mathbb P)}(1)\cdot \delta_{k}^{0}\,\delta_{\MM}^{\mathbb{S}^2}\,,\label{ReductionPrinciples-eq} \\
&c_2(u) = \frac12 H_1(u)\phi(u), \quad c_1(u) = \frac{1}{2}\sqrt{\frac{\pi }{8}}H_1(u)^2\phi (u),\quad c_0(u) = \frac{1}{2} H_1(u)H_2(u)\phi (u) \frac{1}{2\pi },\notag
\end{flalign}
for every $k=0,1,2$ and $u\in \mathbb R$.
While for the excursion area $\LL_2^{f_n}(\EE_u(f_n,\MM))$ the second chaos is immediately proportional to the integral of $H_2(f_n(x))$, on the contrary, more computations are needed in order to show that for the boundary lengths and the Euler-Poincar\'e characteristics it is also the case. Indeed, as one can see in \eqref{2ndchaos-EPC} and \eqref{2ndchaos-BL}, simply using chaotic decomposition, the boundary length depends on both the level and its gradient, while the Euler-Poincar\'e characteristic depends on both first and second derivatives of the field. For the boundary length, the Green's formula (IBP) is used to prove that its second chaos is proportional to the integral of $H_2$ of the level, see  \cite[Section 7.3]{Ro:15b} for the computations on the sphere and both  [remark 2.4]\cite{MPRW:16} and  [Proposition 3.2]\cite{CMR:20} for statements on the torus. A unified discussion can also be found in  \cite{Ro:19}, in particular Section 4, and we show it again in Section \ref{red-BL} for sake of completeness. 
For the EPC, via some analytic computations, \cite{CM:18} and  \cite{CMR:20} show that formula \eqref{ReductionPrinciples-eq} holds. Here we will show that this can also be done using only IBP as for the boundary length, see Section \ref{red-EPC}.

\begin{oss}\label{rem-NA-3}
Note that in formula \eqref{ReductionPrinciples-eq}  the high-energy regime is present only for $k=0$ and $\MM=\S^2$. This means that in all other cases the formula is exact, in the sense that is non-asymptotic. 
\end{oss}

Regarding the nodal case, that is for $u=0$, only the boundary length is dominated by a single chaotic component, the fourth one; this is shown in  \cite{MRW:20} for the sphere and the reduction principle consists in having the dominant term asymptotically proportional to the sample trispectrum, which is the integral of $H_4(f_n)$, $H_4$ being the fourth Hermite polynomial. Analogous results on the torus are  shown  in  \cite{CMR:20}, where, however, the dominant fourth chaos is not proportional to the sample trispectrum.

\subsection{The boundary length}\label{red-BL}

In this section we show the crucial role of Green's formula in order to see the cancellation of the second order Wiener chaos for the boundary length $\LL_1^{f_n}(\EE_u(f_n,\MM))$ when $u=0$. This was shown for the first time in \cite{Ro:15b} in the case of $\MM=\S^2$ and  since the proof is completely independent of the manifold, as pointed out in \cite[Section 4.1.2]{Ro:19}, as well as very short and interesting, we represent it here in a unified way, using our generic notation. In fact, this was the very first reduction principle for LKCs and it was proven via some integration by parts formula. Recalling the Green identity on manifolds
$$
\int_{\MM} f_{n}(x) \Delta f_n(x) dx =  -\int_{\MM} \langle \nabla f_n(x), \nabla f_n(x) \rangle \,dx\,,
$$
one can apply it as follows,
\begin{flalign*}
 \operatorname{Proj}[ \LL_1^{f_n}(\EE_u(f_n,\MM)) |2]&= \sqrt{\frac{\lambda_n}{2}}\Big (\frac{\beta_2 \alpha_{00}}{2!}\int_{\MM} H_2\(f_n(x)\)\,dx\\
 &+ \frac{\beta_0(u) \alpha_{20}}{2!}\int_{\MM} \{H_2\(\wt\partial_1f_n(x)\)+H_2\(\wt\partial_2f_n(x)\)\} dx\\
 &= \sqrt{\lambda_n}\Big ( \frac{\beta_2 \alpha_{00}}{2!}\int_{\mathbb S^2} (f_n(x)^2 - 1)\,dx \\
&+ \frac{\beta_0(u) \alpha_{20}}{2!}\int_{\mathbb S^2} \Big ( \frac{2}{\ell(\ell+1)}\langle \nabla f_n(x)), \nabla f_n(x)\rangle - 2\Big )\,dx \Big )\\
&=\sqrt{\lambda_n}\Big ( \frac{\beta_2 \alpha_{00}}{2!}\int_{\mathbb S^2} (f_n(x)^2 - 1)\,dx \\
& +\frac{\beta_0(u) \alpha_{20}}{2!}\int_{\mathbb S^2} \Big ( - \frac{2}{\ell(\ell+1)}f_n(x) \Delta f_n(x) - 2\Big )\,dx \Big )\\
&=\sqrt{\lambda_n}\Big ( \frac{\beta_2 \alpha_{00}}{2!}\int_{\mathbb S^2} (f_n(x)^2 - 1)\,dx \\
& + \frac{\beta_0(u) \alpha_{20}}{2!}\int_{\mathbb S^2} \Big ( 2 f_n(x)^2  - 2\Big )\,dx \Big )\\
&= \sqrt{\lambda_n}\Big ( \frac{\beta_2 \alpha_{00}}{2!}+  \beta_0(u) \alpha_{20}\Big )\int_{\mathbb S^2} H_2\(f_n(x)\)dx\,,
\end{flalign*}
which is \eqref{ReductionPrinciples-eq} in the case of $k=1$.

\subsection{The Euler-Poincar\'e characteristic}\label{red-EPC}

In this section we want to give an alternative proof of the reduction principles for Euler-Poincar\'e characteristics given in [Theorem 1]\cite{CM:18} and [Theorem 2.4]\cite{CMR:20}. This alternative proof is independent of the manifold $\MM$, except for the constants involved in the computations, that are different on the torus and on the sphere, see also Proposition \ref{prop-chaos-EPC} in Section \ref{LKC and chaos}. For this reason, we will prove \eqref{ReductionPrinciples-eq} in the case $k=0$ setting $\MM=\T^2$, avoiding  to repeat analogous computations but with different constants for $\MM=\S^2$.

\begin{oss}
To be precise, the main difference regarding the computations to reach \eqref{ReductionPrinciples-eq} in the case $k=0$ on the two different manifolds, is the fact that on the sphere one has to consider covariant derivatives instead of \emph{flat} derivatives. The flat geometry of the two-dimensional standard flat torus gives a neater expression for the projection onto the second order Wiener chaos of the Euler-Poincar\'e characteristic, by neater we mean an expression without a reminder -- see \eqref{ReductionPrinciples-eq}.
\end{oss}

From [Section 6.1]\cite{CMR:20}, we know that the projection of $\LL_0^{f_n}(\EE_u(f_n,\T^2))$ onto the second order Wiener chaos can be compactly written as follows, see also [Section 3.2]\cite{CM:18} for the case of $\S^2$, 
\begin{equation}\label{simple-EPC}
\text{\rm Proj}[\mathcal{L}_{0}(n;u)|2]=h_{35}(u;n)\int_{\T^2}Y_{3}(x)Y_{5}(x)dx+\frac{1}{2}\sum_{i=1}^{5}h_{i}(u;n)\int_{\T^2}H_{2}(Y_{i}(x))dx\,,
\end{equation}
where
\begin{flalign*}
&Y_1(x)=\frac{1}{\kappa_1}\partial_{1}f_{n}(x)  \qquad  Y_{2}(x)=\frac{1}{\kappa_1}\partial _{2}f_{n}(x) \qquad Y_{3}(x)=\frac{1}{\kappa_3}\partial _{11}f_{n}(x)\\
&Y_4(x)=\frac{1}{\kappa_4}\partial _{12}f_{n}(x) \qquad Y_5(x)=\frac{1}{\kappa_5}\partial _{22}f_{n}(x)-\frac{\kappa_2}{\kappa_3\kappa_5}\partial _{11}f_{n}(x)\,,
\end{flalign*}
recalling that (see Proposition \ref{prop-chaos-EPC})
\begin{flalign*}
&\kappa_1=\sqrt{\frac{\lambda_n}{2}} \qquad \kappa_2=\frac{\lambda_n}{2\sqrt2}\frac{1-\hat{\mu}_n(4)}{\sqrt{3+\hat{\mu}_n(4)}} \qquad \kappa_3=\frac{\lambda_n}{2\sqrt2}\sqrt{3+\hat{\mu}_n(4)}\\
& \kappa_4=\frac{\lambda_n}{2\sqrt2}\sqrt{1-\hat{\mu}_n(4)} \qquad \kappa_5=\lambda_n\frac{\sqrt{1+\hat{\mu}_n(4)}}{\sqrt{3+\hat{\mu}_n(4)}}\,,
\end{flalign*}
while
\begin{eqnarray*}
h_{35}(u; n )&=& \frac{\lambda_n}{2 \sqrt 2 \pi} \sqrt{1+\hat{\mu}_n(4)} \frac{ u
\phi (u)(1+u^2) + (3+\hat{\mu}_n(4)) \Phi (-u) }{3+\hat{\mu}_n(4)}\,,
\end{eqnarray*}
and
\begin{eqnarray*}
h_{1}(u;n )&=& h_{2}(u;n ) = - \frac{\lambda_n}{4 \pi} u \,\phi(u), \\
h_{3}(u; n ) &=& \frac{\lambda_n}{4 \pi} \left[ \frac{2 u (1+u^2) \phi(u)}{3+\hat{\mu}_n(4) } + \Phi(-u)(1 - \hat{\mu}_n(4)) \right], \\
h_{4}(u;n ) &=& - \frac{\lambda_n}{4 \pi} (1-\hat{\mu}_n(4)) \Phi (-u), \\
h_{5}(u;n) &=& \frac{\lambda_n}{4 \pi} \frac{ u (1+u^2) (1+\hat{\mu}_n(4))\phi (u)} {3+\hat{\mu}_n(4)}.
\end{eqnarray*}

Let us now show that, starting from \eqref{simple-EPC},  we  can simply use integration by  parts to arrive at \eqref{ReductionPrinciples-eq}, also for the case of $k=0$,  that is for the Euler-Poincar\'e characteristic:
\begin{flalign*}
&\text{\rm Proj}[\mathcal{L}_{0}(n;u)|2]=h_{35}(u;n)\int_{\T^2}\frac{1}{\kappa_3}\partial _{11}f_{n}(x)\(\frac{1}{\kappa_5}\partial _{22}f_{n}(x)-\frac{\kappa_2}{\kappa_3\kappa_5}\partial _{11}f_{n}(x)\)dx\\
&+\frac{h_{1}(u;n )}{2}\int_{\T^2}H_{2}\left( \frac{1}{\kappa_1}\partial_{1}f_{n}(x)\right) dx+\frac{h_{1}(u;n )}{2}\int_{%
\T^2}H_{2}\left( \frac{1}{\kappa_1}\partial _{2}f_{n}(x)\right) dx \\
&+\frac{h_{3}(u;n )}{2}\int_{\T^2}H_{2}\left( \frac{1}{\kappa_3}\partial _{11}f_{n }(x)\right) dx+\frac{\kappa_4(u;n )}{2}\int_{%
\T^2}H_{2}\left( \frac{1}{\kappa_4}\partial _{12}f_{n }(x)\right) dx \\
&+\frac{\kappa_5(u;n )}{2}\int_{\T^2}H_{2}\left( \frac{1}{\kappa_5}\partial _{22}f_{n }(x)-\frac{\kappa_2}{\kappa_3\kappa_5}\partial _{11}f_{n}(x)\right) dx \\
&=h_{35}(u;n )\int_{\T^2}\left( \frac{1}{\kappa_3\kappa_5}\partial _{22}f_{n}(x)\partial _{11}f_{n}(x)-\frac{\kappa_2}{\kappa_3^2\kappa_5}\left( \partial _{11}f_{n}(x)\right) ^{2}\right) \,dx \\
&+\frac{\kappa_1(u;n )}{2}\int_{\T^2}\left( \frac{1}{k
_{1}^{2}}\left(\partial _{1}f_{n}(x)\right) ^{2}-1\right) dx+\frac{\kappa_1(u;n )}{2}\int_{\T^2}\left( \frac{1}{\kappa_1^{2}}%
\left( \partial _{2}f_{n }(x)\right) ^{2}-1\right) dx \\
&+\frac{\kappa_3(u;n )}{2}\int_{\T^2}\left( \frac{1}{k
_{3}^{2}}\left( \partial _{11}f_{n }(x)\right) ^{2}-1\right) dx+
\frac{\kappa_4(u;n )}{2}\int_{\T^2}\left( \frac{1}{\kappa_4^{2}
}\left( \partial _{12}f_{n }(x)\right) ^{2}-1\right) dx \\
&+\frac{\kappa_5(u;n )}{2}\int_{\T^2}\left( \frac{1}{k
_{5}^{2}}\left( \partial _{22}f_{n }(x)-\frac{\kappa_2}{k
_{3}}\partial _{11}f_{n }(x)\right) ^{2}-1\right) dx\,.
\end{flalign*}
Now we use the fact that 
$$
\int_{\T^2}\left( \partial_{12}f_{n}(x)\right) ^{2}\,dx=\int_{\T^2}\partial_{11}f_{n}(x)\partial_{22}f_{n}(x)\,dx
$$
to have 
\begin{flalign*}
&\mathtt{Proj}[\chi (A_{u}(f_{n};\T^2))|2] \\
&=h_{35}(u;n )\int_{\T^2}\left( \frac{1}{\kappa_3k
_{5}}\partial_{22}f_{n}(x)\partial_{11}f_{n}(x)-\frac{
\kappa_2}{\kappa_3^{2}\kappa_5}\left( \partial_{11}f_{n}(x)\right) ^{2}\right) \,dx \\
&+\frac{h_{1}(u;n )}{2}\int_{\T^2}\left( \frac{1}{k
_{1}^{2}}\left( \partial_1f_{n}(x)\right) ^{2}-1\right) dx+\frac{
h_{1}(u;n )}{2}\int_{\T^2}\left( \frac{1}{\kappa_1^{2}}
\left( \partial_2f_{n}(x)\right) ^{2}-1\right) dx \\
&+\frac{h_{3}(u;n )}{2}\int_{\T^2}\left( \frac{1}{k
_{3}^{2}}\left( \partial_{11}f_{n}(x)\right) ^{2}-1\right) dx+
\frac{h_{4}(u;n )}{2}\int_{\T^2}\left( \frac{1}{\kappa_4^{2}
}\partial_{11}f_{n}(x)\partial_{22}f_{n}(x)-1\right) dx
\\
&+\frac{h_{5}(u;n )}{2}\int_{\T^2}\left( \frac{1}{k
_{5}^{2}}\left( \partial_{22}f_{n}(x)-\frac{\kappa_2}{k
_{3}}\partial_{11}f_{n}(x)\right) ^{2}-1\right) dx \\
&=\left[ \frac{1}{\kappa_3\kappa_5}h_{35}(u;n )+\frac{1}{k
_{4}^{2}}\frac{h_{4}(u;n )}{2}\right] \int_{\T^2}\partial
_{2,x}^{2}f_{n}(x)\partial_{11}f_{n}(x)dx \\
&+\left[ \frac{h_{3}(u;n )}{2\kappa_3^{2}}-\frac{k
_{2}\,h_{35}(u;n )}{\kappa_3^{2}\kappa_5}\right] \int_{\T^{2}}\left( \partial_{11}f_{n}(x)\right) ^{2}\,dx+\frac{%
h_{1}(u;n )}{2\kappa_1^{2}}\int_{\T^2}\Vert \nabla f_{n}(x)\Vert ^{2}\,dx \\
&+\frac{h_{5}(u;n )}{2}\int_{\T^2}\left[ \frac{1}{k
_{5}^{2}}\left( \partial_{22}f_{n}(x)\right) ^{2}+\frac{k
_{2}^{2}}{\kappa_3^{2}\kappa_5^{2}}\left( \partial_{11}f_{n}(x)\right) ^{2}-2\frac{\kappa_2}{\kappa_3\kappa_5^{2}}%
\,\partial_{22}f_{n}(x)\partial_{11}f_{n}(x)\right] dx \\
&- \left[ h_{1}(u;n )+\frac{h_{3}(u;n )+h_{4}(u;n
)+h_{5}(u;n )}{2}\right]  \\
&=A(u;n) \int_{\T^2}\partial
_{2,x}^{2}f_{n}(x)\partial_{11}f_{n}(x)dx+B(u;n) \int_{
\T^2}\left( \partial_{11}f_{n}(x)\right) ^{2}\,dx \\
&+C(u;n)\int_{\T^2}\Vert \nabla
f_{n}(x)\Vert ^{2}\,dx+D(u;n)
\int_{\T^2}\left( \partial_{22}f_{n}(x)\right) ^{2}dx - E(u;n) \,.
\end{flalign*}
After using integration by parts, now we just have to compute the constants in front of the integral terms:
\begin{flalign*}
A(u;n)&=\frac{1}{\kappa_3\kappa_5}h_{35}(u;n )+\frac{1}{\kappa_4^{2}}\frac{h_{4}(u;n )}{2}-\frac{\kappa_2}{\kappa_3\kappa_5^{2}}h_{5}(u;n )=\\
&=\frac{2\sqrt2}{\lambda_n^2\sqrt{1+\hat{\mu}_n(4)}}\frac{\lambda_n}{2 \sqrt 2 \pi} \sqrt{1+\hat{\mu}_n(4)} \frac{ u
\phi (u)(1+u^2) + (3+\hat{\mu}_n(4)) \Phi (-u) }{3+\hat{\mu}_n(4)}\\
&-\frac{8}{\lambda_n^2(1-\hat{\mu}_n(4))}\frac{\lambda_n}{8 \pi} (1-\hat{\mu}_n(4)) \Phi (-u)\\
&-\frac{\lambda_n}{2\sqrt2}\frac{1-\hat{\mu}_n(4)}{\sqrt{3+\hat{\mu}_n(4)}}\frac{2\sqrt2}{\sqrt{3+\hat{\mu}_n(4)}{\lambda_n}}
\frac{3+\hat{\mu}_n(4)}{\lambda_n^2(1+\hat{\mu}_n(4))}\frac{\lambda_n}{4 \pi} \frac{ u (1+u^2) (1+\hat{\mu}_n(4))\phi (u)} {3+\hat{\mu}_n(4)}\\
&=\frac{1}{\lambda_n\pi} \frac{u\phi (u)(1+u^2) + (3+\hat{\mu}_n(4)) \Phi (-u) }{3+\hat{\mu}_n(4)}-\frac{1}{\lambda_n\pi}\Phi (-u)\\
&-\frac{1-\hat{\mu}_n(4)}{3+\hat{\mu}_n(4)}\frac{1}{4 \lambda_n \pi}  u (1+u^2) \phi (u)=\frac{u\phi (u)(1+u^2)}{4\lambda_n\pi} 
\end{flalign*}

\begin{flalign*}
B(u;n)&=\frac{h_{3}(u;n )}{2\kappa_3^{2}}-\frac{\kappa_2\,h_{35}(u;n )}{\kappa_3^{2}\kappa_5}+\frac{\kappa_2^{2}}{\kappa_3^{2}\kappa_5^{2}}\frac{h_{5}(u;n )}{2}\\
&=\frac{1}{\kappa_3^{2}} \times \(\frac{h_{3}(u;n )}{2}-\frac{\kappa_2\,h_{35}(u;n )}{\kappa_5}+\frac{\kappa_2^{2}}{\kappa_5^{2}}\frac{h_{5}(u;n )}{2}\)\\
&=\frac{8}{\lambda_n^2(3+\hat{\mu}_n(4))}\times\\
&\(\frac{\lambda_n}{8 \pi} \left[ \frac{2 u (1+u^2) \phi(u)}{3+\hat{
\mu}_n(4) } + \Phi(-u)(1 - \hat{\mu}_n(4)) \right]\right.\\
&-(1-\hat{\mu}_n(4))\(\frac{\lambda_n}{8 \pi}  \frac{ u
\phi (u)(1+u^2) + (3+\hat{\mu}_n(4)) \Phi (-u) }{3+\hat{\mu}_n(4)}\)\\
&\left.+\frac{(1-\hat{\mu}_n(4))^2}{3+\hat{\mu}_n(4)}\frac{\lambda_n}{64 \pi}  u (1+u^2) \phi (u)\)\\
&=\frac{u (1+u^2) \phi (u)}{8\lambda_n\pi}
\end{flalign*}
\begin{flalign*}
C(u;n)=\frac{1}{\kappa_5^{2}}\frac{h_{5}(u;n)}{2}=\frac{u (1+u^2) \phi (u)}{8 \lambda_n \pi} 
\qquad
D(u;n)=\frac{h_{1}(u;n)}{2\kappa_1^{2}}=-\frac{u \,\phi(u)}{4\pi}
\end{flalign*}
\begin{flalign*}
E(u;n)&=h_{1}(u;n )+\frac{h_{3}(u;n )+h_{4}(u;n)+h_{5}(u;n )}{2}\\
&=- \frac{\lambda_n}{4 \pi} u \,\phi(u)+\frac{\lambda_n}{8 \pi} \left[ \frac{2 u (1+u^2) \phi(u)}{3+\hat{
\mu}_n(4) } \right]+\frac{\lambda_n}{8 \pi}  u (1+u^2)\phi (u)-\frac{\lambda_n}{8 \pi}\frac{2u(1+u^2)\phi(u)}{3+\hat{\mu}_n(4)}\\
&=\frac{\lambda_n}{8 \pi}  u (u^2-1)\phi (u)=\frac{\lambda_n}{8 \pi}  H_1(u) H_2(u)\phi (u)
\end{flalign*}
As a consequence, we easily have that
\begin{flalign*}
\mathtt{Proj}[\chi (A_{u}(f_{n};\T^2))|2]&=\frac{u\phi (u)(1+u^2)}{4\lambda_n\pi} \int_{\T^2}\partial
_{2,x}^{2}f_{n}(x)\partial_{11}f_{n}(x)dx \\
&+\frac{u (1+u^2) \phi (u)}{8 \lambda_n \pi}  \int_{
\T^2}\left( \partial_{11}f_{n}(x)\right) ^{2}\,dx \\
&-\frac{u \,\phi(u)}{4\pi}\int_{\T^2}\Vert \nabla
f_{n}(x)\Vert ^{2}\,dx+\frac{u (1+u^2) \phi (u)}{8 \lambda_n \pi} 
\int_{\T^2}\left( \partial_{22}f_{n}(x)\right) ^{2}dx \\
&- \frac{\lambda_n}{8 \pi}  H_1(u) H_2(u)\phi (u) \\
&=\frac{u (1+u^2) \phi (u)}{8 \lambda_n \pi}\int_{\T^2} \(\Delta_{\T^2}f_{n}(x)\)^2dx-\frac{u \,\phi(u)}{4\pi}\int_{\T^2}\Vert \nabla
f_{n}(x)\Vert ^{2}\,dx\\
&-\frac{\lambda_n}{8 \pi}  H_1(u) H_2(u)\phi (u)
\end{flalign*}
and recalling the basic (Green-Stokes) identity
$$
\int_{\T^2}\norm{\nabla f_n}^{2}\,dx=-\int_{\T^2} f_n \Delta_{\T^2}f_n\,dx\,
$$
we obtain
\begin{flalign*}
\mathtt{Proj}[\chi (A_{u}(f_{n};\T^2))|2]
&=\frac{u (1+u^2) \phi (u)}{8 \lambda_n \pi}\int_{\T^2} \(\lambda_n f_{n}(x)\)^2dx\\
&+\frac{u \,\phi(u)}{4\pi}\int_{\T^2}
f_{n}(x)\Delta_{\T^2}f_{n}(x)\,dx-\frac{\lambda_n}{8 \pi}  H_1(u) H_2(u)\phi (u)\\
&=\frac{u (1+u^2) \phi (u)}{8 \lambda_n \pi}\int_{\T^2} \lambda_n^2 f_{n}(x)^2dx-\frac{2u \,\phi(u)}{8\pi}\int_{\T^2}
\lambda_n f_{n}(x)^2\,dx\\
&-\frac{\lambda_n}{8 \pi}  H_1(u) H_2(u)\phi (u)\\
&=\frac{H_1(u) H_2(u) \phi (u)\lambda_n}{8\pi}\int_{\T^2}H_2\(f_{n}(x)\)dx\,,
\end{flalign*}
which is the desired formula.

 \bibliographystyle{amsplain}
 \bibliography{Bibliography}

\end{document}